\documentclass[10pt,twoside]{article}
\usepackage{amsthm,amsbsy,amsfonts,amssymb,amsmath,amscd}
\usepackage{latexsym}
\usepackage{fontenc}
\usepackage[mathscr]{euscript}

\usepackage[dvips]{epsfig}
\usepackage{Latex-document}

\DeclareMathAlphabet{\scrb}{U}{eus}{b}{n}
\theoremstyle{plain}

\newtheorem{Thm}{Theorem}[section]

\newtheorem{Prop}{Proposition}[section]

\theoremstyle{remark}

\newtheorem{Thme}[Prop]{\bf Theorem}

\def\bcw{\mathbin{\bigcirc\mkern-15mu\wedge}}

\renewcommand{\thesection}{\arabic{section}}
\renewcommand{\thesubsection}{\thesection.\arabic{subsection}.}

\def\volumeno{I}

\title{\bf Non-linear Partial Differential Equations \vskip -2mm in Conformal
Geometry\thanks{Research of Chang is supported in part by NSF Grant
DMS-0070542. Research of Yang is supported in part by NSF Grant
DMS-0070526.}\vskip 6mm}

\author{{\bf Sun-Yung Alice Chang}\thanks{Department of Mathematics,
Princeton University, Princeton, NJ 08544, USA. E-mail:
chang@math.princeton.edu} \quad Paul C. Yang\thanks{Department of
Mathematics, Princeton University, Princeton, NJ 08544, USA.
E-mail: yang@math.princeton.edu}\vspace*{-0.5cm}}

\date{\vspace{-8mm}}

\begin{document}

\maketitle \thispagestyle{first} \setcounter{page}{189}

\setcounter{section}{-1}

\section{\hskip -1em . \hskip 0.3em Introduction} \setzero
\vskip-5mm \hspace{5mm}

In the study of conformal geometry, the method of elliptic partial
differential equations is playing an increasingly significant role.
Since the solution of the Yamabe problem, a family of conformally
covariant operators (for definition, see section 2) generalizing the
conformal Laplacian, and their associated conformal invariants have been
introduced. The conformally covariant powers of the Laplacian form a
family $P_{2k}$ with $k \in \mathbb N$ and $k \leq \frac{n}{2}$ if the
dimension $n$ is even. Each $P_{2k}$ has leading order term $(-\, \Delta
)^k$ and is equal to $ (- \, \Delta) ^k$ if the metric is flat.

The curvature equations associated with these $P_{2k}$ operators
are of interest in themselves since they exhibit a large
group of symmetries. The analysis of these equations is of necessity more
complicated,
it typically requires the derivation of an
optimal Sobolev or Moser-Trudinger inequality that always occur at a critical
exponent. A common feature is the presence of blowup or bubbling
associated to the noncompactness of the conformal group. A number of
techniques have been introduced to study the nature of blowup, resulting
in a well developed technique to count the topological degree of such
equations.

The curvature invariants (called the Q-curvature) associated to such operators are
also of higher order. However, some of the invariants are closely related
with the Gauss-Bonnet-Chern integrand in even dimensions, hence of intrinsic
interest to geometry. For example, in dimension four, the
finiteness of the $Q$-curvature integral can be used to conclude
finiteness of topology.
In addition, the symmetric functions of the Ricci tensor appear in natural fashion
as the lowest order terms of these curvature invariants, these equations offer
the possibility to analyze the Ricci tensor itself. In particular, in dimension
four the sign of the $Q$-curvature integral can be used to conclude the
sign of the Ricci tensor. Therefore there is ample
motivation for the study of such equations.

In the following sections we will survey some of the development
in the area that we have been involved. We gratefully acknowledge
the collaborators that we were fortunate to be associated with.

\section{\hskip -1em . \hskip 0.3em Prescribing Gaussian curvature on compact
surfaces and the Yamabe problem} \setzero \vskip-5mm \hspace{5mm}

In this section we will describe some second order elliptic equations
which have played important roles in conformal geometry.

On a compact surface $(M, g)$ with a Riemannian metric $g$, a natural
curvature invariant associated with the Laplace operator $\Delta = \Delta
_g$ is the Gaussian curvature $K=K_g$. Under the conformal change of
metric $g_w = e^{2w} g$, we have
\begin{equation}
- \Delta w \, + \, K = K_w e^{2w} \,\, on \,\, M
\end{equation}
where $K_w$ denotes the Gaussian curvature of $(M, g_w)$. The classical
uniformization theorem to classify compact closed surfaces can be viewed
as finding solution of equation (1.1) with $K_w \equiv -1$, $0$, or
$1$
according to the sign of $\int K dv_g$. Recall that the
Gauss-Bonnet theorem states
\begin{equation}
\int_M  K_w \, dv_{g_w} = 2 \pi \, \chi(M)
\end{equation}
where $\chi(M)$ is the Euler characteristic of $M$, a topological
invariant. The variational functional with
(1.1) as Euler equation for $K_w = constant $  is thus given by
\begin{equation}
J[w] = \int_M  |\nabla w |^2 dv_g + 2 \int_M  K w dv_g -(\int_M K
dv_g) \log \frac { \int_M dv_{g_w} }{\int_M dv_g}.
\end{equation}

When the surface $(M, g)$ is the standard 2-sphere $S^2$ with the
standard canonical metric, the problem of prescribing Gaussian
curvature on $S^2$ is commonly known as the Nirenberg problem. For
general compact surface $M$, Kazdan and Warner (\cite {KW}) gave a
necessary and sufficient condition for the function when $\chi (M)
= 0 $ and some necessary condition for the function when $\chi (M)
< \, 0$. They also pointed out that in the case when $\chi (M) >
0$, i.e. when $ (M, g) = (S^2,g_c )$, the standard 2-sphere with
the canonical metric $ g= g_c$, there is an obstruction for the
problem:
\begin{equation}
\int _{S^2} \nabla K_{w} \cdot \nabla x \,\, e^{2w} dv_g  = 0
\end{equation}
where $x$ is any of the ambient coordinate function. Moser  (\cite{M-2})
realized that this implicit integrability condition is satisfied if the conformal factor
has antipodal symmetry. He proved for an even function $f$, the  only
necessary condition for (1.1) to be
solvable with $K_w = f$ is that $f$ be positive somewhere.
An important tool introduced by Moser
is the following inequality (\cite {M-1}) which is a sharp form of an earlier result of
Trudinger (\cite {Tr}) for the limiting Sobolev embedding of $W_0^{1,2}$ into the
Orlicz space $e^{L^2}$:
Let $w$ be a smooth function on the $2$-sphere satisfying the
normalizing conditions:
$\int _{S^2}|\nabla w|^2 dv_g \leq 1$ and $\bar w=0$ where
$\bar w$ denotes the mean value of $w$,
then
\begin{equation}
\int _{S^2}e^{\beta w^2}dv_g \leq C
\end{equation}
where $\beta \leq 4\pi$ and $C$ is a fixed constant and $4 \pi$ is
the best constant. If $w$ has antipodal symmetry
then the inequality holds for $\beta \leq 8\pi$.

Moser has also established a similar inequality
for functions $u$ with compact support on bounded domains in the
Euclidean space
$\mathbb R^n$ with the $W^{1,n}$ energy norm  $\int |\nabla u|^n dx$ finite.
Subsequently, Carleson and Chang (\cite{CC}) found that, contrary to the situation
for Sobolev embedding, there is an extremal function realizing
the maximum value of the inequality of Moser when the domain is the unit ball
in Euclidean space.   This fact
remains true for simply connected domains in the plane
 (Fl\"ucher \cite{Fl}),
and for some
domains in the n-sphere  (Soong \cite{So}).

Based on the inequality of Moser and subsequent work of Aubin
(\cite {Au-2}
and
Onofri (\cite {On}), we devised a degree
count (\cite {CY-1}, \cite{CY-2}, \cite{CGY-1}) associated to the function
$f$ and the Mobius group on the $2$-sphere, that is motivated
by the Kazdan-Warner condition (1.4).
This degree actually computes  the  Leray-Schauder
degree of the equation (1.1) as a nonlinear Fredholm equation.
In the special case that $f$ is a Morse function satisfying the condition
$\Delta f(x)\neq 0$ at the critical points $x$ of $f$, this degree can be expressed as:
\begin{equation}
\sum _{\nabla f(q)=0,\Delta f(q)<0}(-1)^{ind (q)} -1.
\end{equation}
The latter degree count is also obtained later by Chang-Liu (\cite{ChL})
and Han (\cite{H}).

There is another interesting geometric interpretation
of the functional $J$ given by Ray-Singer (\cite{RS}) and Polyakov
(\cite{Po}); (see also Okikiolu \cite {Ok-1})
\begin{equation}
J[w]=  12 \pi \,\,\, \log \,\,(\frac{\det \,\, \Delta _g}{\det
\,\, \Delta _{g_w}})
\end{equation}
for metrics $g_w$ with the volume of $g_w$ equals the volume of
$g$; where the determinant of the Laplacian $det \,\, \Delta _g$
is defined by Ray-Singer via the ``regularized'' zeta function. In
\cite{On}, (see also Hong \cite{Ho}), Onofri established the sharp
inequality that on the 2-sphere $J [w] \geq 0$  and $ J [w] = 0 $
precisely for conformal factors $w$ of the form $e^{2w}g_0=
T^*g_0$ where $T$ is a Mobius transformation of the 2-sphere.
Later Osgood-Phillips-Sarnak (\cite{OPS-1}, \cite{OPS-2}) arrived
at the same sharp inequality in their study of heights of the
Laplacian. This inequality also plays an important role in their
proof of the $C^{\infty}$ compactness of isospectral metrics on
compact surfaces.

The formula of Polyakov-Ray-Singer has been generalized to
manifolds of dimension greater than two in many different
settings; one of which we will discuss in section 2 below. There
is also a general study of extremal metrics for $ det \,\,
\Delta_g $ or $ det \,\, L_g$ for metrics $g$  in the same
conformal class with a fixed volume or for all metrics with a
fixed volume(\cite {Be}, \cite {BCY}, \cite {Br-3}, \cite {R},
\cite {Ok-2}). A special case of the remarkable results of
Okikiolu (\cite {Ok-2}) is that among all metrics with the same
volume as the standard metric on
 the 3-sphere, the standard canonical metric is a local maximum for the
functional $det \,\,\Delta_g$.

More recently, there is an extensive study of a generalization
of the equation
(1.1)
to compact Riemann surfaces. Since Moser's argument is readily applicable to a
compact surface $(M,g)$, a lower bound for similarly defined
functional $J$ on $(M, g)$ continues to hold in that situation.
The Chern-Simons-Higgs equation in the Abelian case is
given by:
\begin{equation}
\Delta w =  \rho e^{2w}(e^{2w} -1)+ 2 \pi \sum_{i=1}^N
\delta_{p_i}.
\end{equation}
A closely related equation is the  mean field equation:
\begin{equation}
\Delta w+ \rho (\frac{he^{2w}}{\int he^{2w}} -1)=0,
\end{equation}
where $\rho$ is a real parameter that is allowed to vary.

There is active development on these equations by several group of
researchers including (\cite{CaYa}, \cite{DJLW}, \cite{T}, \cite{ST},
\cite {CL-1}).
%In particular in (\cite{CL-1})
%C.-C. Chen and C.-S. Lin gave a degree
%count for the mean field equation for all parameter
%values of $\rho$ on all
%compact surfaces.

On manifolds $(M^n, g)$ for n greater than two,  the conformal
Laplacian $L_g$ is defined as
$L_g = - c_n \Delta_g + R_g$ where $c_n = \frac {4 (n-1)}{n-2}$, and $R_g$ denotes the scalar
curvature of the metric $g$. An analogue of equation (1.1) is the equation, commonly referred to as the Yamabe equation, which relates
the scalar
curvature under conformal change of metric to the background metric. In this case, it is
convenient to denote the conformal metric as $ \bar g = u ^{ \frac {4}{n-2}} g$ for some positive function $u$, then the equation becomes
\begin{equation}
L_g u \,\, = \,\, \bar R \,\, u^{\frac {n+2}{n-2}}.
\end{equation}
The famous Yamabe problem to solve (1.10) with $  \bar R $ a constant has been
 settled by Yamabe (\cite {Y}), Trudinger (\cite {Tr-2}), Aubin (\cite {Au-1}) and
Schoen (\cite {S}). The corresponding problem to prescribe scalar
curvature has been intensively studied in the past decades by
different groups of mathematicians, we will not be able to survey
all the results here. We will just mention that the degree theory
for existence of solutions on the $n$-sphere has been achieved by
Bahri-Coron (\cite {BC}), Chang-Gursky-Yang (\cite {CGY-1}) and
Schoen-Zhang (\cite {SZ}) for $n=3$ and under further constraints
on the functions for $n \geq 4$ by Y. Li (\cite {L-2}) and by C-C.
Chen and C.-S. Lin (\cite {CL-2}).

\def\intl{\int\kern-9pt\hbox{$\backslach$}}

\section{\hskip -1em . \hskip 0.3em Conformally covariant differential
operators and the $Q$-curvatures} \setzero
\vskip-5mm \hspace{5mm}

It is well known that in dimension two, under the conformal change of
metrics $g_w= e^{2w}g$, the associated Laplacians are related by
\begin{equation}
\Delta_{g_w}= e^{-2w}\Delta _g.
\end{equation}
Similarly on $(M^n, g)$, the conformal Laplacian $L= -\frac{4(n-1)}{n-2}\Delta + R$
transforms under the conformal change of metric $\bar g=u^{\frac{4}{n-2}}g$:
\begin{equation}
L_{\bar g}= u^{-\frac{n+2}{n-2}}L_g(u\cdot)
\end{equation}

 In general, we call a metrically defined operator $A$ conformally covariant
 of bidegree $(a, b)$,  if under the conformal change of metric
 $g_\omega = e^{2\omega }g$, the pair of corresponding operators
 $A_\omega $ and $A$ are related by
$$
A_\omega (\varphi ) = e^{-b\omega } A(e^{a\omega }\varphi )\quad
\text{\rm for all}\quad \varphi \in C^\infty (M^n) \, .
$$
Note that in this notation, the conformal Laplacian opertor is
conformally covariant of bidegree $(\frac{n-2}{2}, \frac{n+2}{2})$.

There are many operators besides the Laplacian $\Delta $
 on compact surfaces and the conformal Laplacian $L$ on general compact
 manifold of dimension greater than two  which have the conformal covariance
 property.  We begin with the fourth order operator on
 $4$-manifolds discovered by Paneitz (\cite {Pa}) in 1983 (see also
\cite {ES}):
 $$
 P\varphi \equiv \Delta ^2 \varphi + \delta
 \left( \frac 23 Rg - 2 \text{\rm Ric}\right) d \varphi
 $$
 where $\delta $ denotes the divergence, $d$ the deRham differential and
 $Ric$ the Ricci tensor of the metric. The {\it Paneitz} operator $P$
 (which we will later denote by $P_4$) is conformally covariant of bidegree
 $(0, 4)$ on 4-manifolds, i.e.
 $$
 P_{g_w} (\varphi ) = e^{-4\omega } P_g (\varphi )\quad
 \text{\rm for all}\quad \varphi \in C^\infty (M^4) \, .
 $$
 More generally, T. Branson (\cite{Br-1}) has extended the definition of the fourth
 order operator to general dimensions $n\neq 2$; which we call the conformal Paneitz operator:
\begin{equation}
P_4^n= \Delta ^2 + \delta \left( a_n Rg + b_n \text{\rm Ric} \right)
 d  + \frac{n-4}{2} Q_4^n
 \end{equation}
 where
 \begin{equation}
 Q_4^n = c_n |Ric |^2 + d_n R^2 - \frac{1}{2(n-1)}\Delta R,
 \end{equation}
 and
$$a_n = \frac {(n-2)^2 + 4}{2(n-1)(n-2)},\;
 b_n= - \frac 4{n-2},\;
 c_n = -\frac{2}{(n-2)^2},\;
 d_n =\frac{n^3-4n^2+16n-16}{8(n-1)^2(n-2)^2}.
$$
 The conformal Paneitz operator is conformally covariant of bidegree $(\frac{n-4}{2},\frac{n+4}{2})$.
 As in the case of the second order conformally covariant operators, the fourth
 order Paneitz operators have associated fourth order curvature invariants $Q$:
 in dimension $n=4$ we write the conformal metric $ g_w = e^{2w}g$;
$ Q = Q_g = \frac {1}{2}(Q_4^4)_g$, then
 \begin{equation}
 Pw + 2Q \, = \, 2 Q_{g_w} e^{4w}
 \end{equation}
 and in dimensions $n \neq 1,2,4$
 we write the conformal metric as $\bar g=u^{\frac{4}{n-4}}g$:
 \begin{equation}
 P_4^n u \, = \, \bar Q_4^n u^{\frac{n+4}{n-4}}.
 \end{equation}
 In dimension $n=4$ the $Q$-curvature equation is closely connected to the Gauss-Bonnet-Chern formula:
 \begin{equation}
 4\pi^2 \chi (M^4) = \int (Q + \frac {1}{8} |W|^2 ) \,\, dv
 \end{equation}
 where $W$ denotes the Weyl tensor, and the quantity $|W|^2dv$
 is a pointwise
 conformal invariant. Therefore the $Q$-curvature integral $\int Q dv$ is a conformal
 invariant.
 The basic existence theory for the $Q$-curvature equation is outlined in \cite {CY-5}:
 \begin{Thm}
 If $\int Q dv < 8 \pi ^2$ and the P operator is positive except for constants,
 then equation (2.5) may be solved with $ Q_{g_w}$ given by a constant.
 \end{Thm}
 It is remarkable that the conditions in this existence theorem are shown by
 M. Gursky (\cite{G-2}) to be
 a consequence of the assumptions that $(M,g)$ has positive Yamabe invariant\footnote{The Yamabe invariant $Y(M, g)$ is defined to be
$ Y(M, g) \equiv \inf_w \frac{\int_M R_{g_w} dv_{g_w}}{{vol(g_w)}^{\frac{n-2}{n}}}$; where $n$ denotes the dimension of M. $Y(M, g)$ is confomally invariant and the sign of $Y (M, g)$ agrees with that of the first eigenvalue of $L_g$.}, and that
 $\int Q dv >0$. In fact, he proves that under these conditions $P$ is a positive
 operator and $\int Q dv \leq 8 \pi ^2$
and that equality can hold only if $(M,g)$ is conformally equivalent to the
standard 4-sphere. This latter fact may be viewed as the analogue of the
positive mass theorem that is the source for the basic compactness result
for the $Q$-curvature equation as well as the associated
fully nonlinear second order equations that we discuss in section
 4.
Gursky's argument is based on a more general existence
result in which we consider a family of 4-th order equations
\begin{equation}
\gamma _1 |W|^2 + \gamma _2 Q - \gamma _3 \Delta R = \bar k \cdot
\text{Vol}^{-1}
\end{equation}
 where
 $\bar k= \int (\gamma _1 |W|^2 + \gamma _2 Q) dv$.
 These equations typically arise as the Euler equation of the functional
 determinants.
 For a conformally covariant operator $A$ of bidegree $(a,b)$ with $b-a=2$
Branson and Orsted (\cite{BO}) gave an explicit
computation
of the normalized form of $\log\,\frac{\det\,A_w}{\det\,A}$ which may
be expressed as:
\begin{equation}
F[w]=\gamma_1I[w]+\gamma_2II[w]+\gamma_3III[w]
\end{equation}
where $\gamma_1,\gamma_2,\gamma_3$ are constants depending only on $A$ and$$
\aligned
I[w]&=4\int\,|W|^2wdv-\left(\int\,|W|^2dv\right)\,\,\log\,
\frac {\int e^{4w}dv} { \int dv }, \\
II[w]&=\langle
Pw,w\rangle+4\int\,Qwdv-\left(\int\,Qdv\right)\,\log\,
\frac {\int  e^{4w}dv}{ \int dv} ,\\
III[w]&= \frac {1}{3} \left ( \int R_{g_w}^2 dv_{g_w} - \int R^2 dv
\right).
\endaligned
$$

In \cite{CY-5}, we gave the general existence result:
\begin{Thm}
If the functional $F$ satisfies
$\gamma_2 >0,\,\,\gamma_3 >0$, and $ \bar k< 8 \gamma_2 \pi^2$, then
$\mathop{\inf}\limits_{w\in W^{2,2}}F[w]$ is attained by some
function $w_d$ and the metric $g_d=e^{2w_d}g_0$ satisfies the
equation
\begin{equation}
\gamma_1\,|W|^2+\gamma_2\,Q_d-\gamma_3\triangle_dR_d \, = \, \bar k\cdot\text
{Vol}(g_d)^{-1}.
\end{equation}
Furthermore, $g_d$ is smooth.
\end{Thm}

This existence result is based on extensions of Moser's inequality
by Adams (\cite {A}, on manifolds \cite {Fo}) to operators of
higher order. In the special case of $(M^4, g)$, the inequality
states that for functions in the Sobolev space $W^{2,2}(M)$ with
\linebreak $\int_M (\Delta w)^2 dv_g \leq 1$, and $\bar w = 0$, we
have
\begin{equation}
\int_M e^{ 32 \pi^2 w^2} dv_g \leq C,
\end{equation}
for some constant $C$. The regularity for minimizing solutions was
first given in \cite{CGY-2}, and later extended to all solutions
by Uhlenbeck and Viaclovsky (\cite {UV}). There are several
applications of these existence result to the study of conformal
structures in dimension $n=4$. In section 4 we will discuss the
use of such fourth order equation as regularization of the more
natural fully nonlinear equation concerned with the Weyl-Schouten
tensor. Here we will mention some elegant application by M. Gursky
(\cite{G-1}) to characterize a number of extremal conformal
structures.

\begin{Thm}
Suppose $(M,g)$ is a compact oriented manifold of dimension four with positive Yamabe invariant.
\item{(i)}  If $\int Q_g dv_g = 0$, and if $M$ admits a non-zero harmonic
\text{\rm 1-}form, then $(M, g)$ is conformal equivalent to a quotient of the
product space $S^3 \times \mathbb R$. In particular $(M, g)$ is locally
conformally flat.
\item{(ii)} If $b_2^+>0$ (i.e. the intersection form has a positive element),
then with respect to the decompostion of the Weyl tensor into
the self dual and anti-self dual components $W= W^+ \oplus W^-$,
\begin{equation}
\int_M |W_g^+|^2 dv_g \geq \frac{4 \pi ^2}{3}(2 \chi + 3\tau),
\end{equation}
where $\tau$ is the signature of $M$. Moreover the equality holds if and only if $g$ is conformal to a
(positive) Kahler-Einstein
metric.
\end{Thm}

In dimensions higher than four, the analogue of the Yamabe equation for the
fourth order Paneitz equation is being investigated by a number
of authors. In particular, Djadli-Hebey-Ledoux (\cite{DHL}) studied the
question of coercivity of the operators $P$ as well as the positivity
of the solution functions,
Djadli-Malchiodi-Ahmedou (\cite{DMA}) have studied the
blowup analysis of the Paneitz equation.
In dimension three, the fourth order Paneitz equation involves a negative exponent, there is now
an existence result (\cite{XY}) in case the Paneitz operator is positive.

 In general dimensions there is an extensive theory of local conformal
 invariants according to the theory of Fefferman and Graham (\cite{FG-1}). For manifolds of general dimension $n$, when $n$ is even, the existence
of a $n$-th order operator $P_n$ conformally covariant of
 bidegree $(0, n)$
 was verified in \cite {GJMS}.  However it is only
 explicitly known on the standard
 Euclidean space $\mathbb R ^n$ and hence on the standard sphere $S^n$.
 For all $n$, on $(S^n, g)$, there also exists an $n$-th order (pseudo)
differential operator $\mathbb P_n$ which is the pull back via
sterographic projection of the operator $(- \Delta)^{n/2}$ from $\mathbb R^n$
with Euclidean metric to $(S^n, g)$.  $\mathbb P_n$ is conformally covariant
of bi-degree (0, n), i.e. $(\mathbb P_n)_w = e^{-nw} \mathbb P_n$.
The explicit formulas for $\mathbb P_n$ on $S^n$ has been computed
in Branson (\cite {Br-3}) and Beckner (\cite {Be}):
\begin{equation}
\begin{cases}
\text{\rm For}\quad n \quad \text{\rm even}\quad\ \Bbb P_n &= \prod
_{k=0}^{\frac {n-2}2}
(- \Delta + k(n-k-1)),\\
\text{\rm For}\quad n \quad \text{\rm odd}\quad\
\Bbb P_n &=
\left(-\Delta + \left(\frac {n-1}2\right)^2\right)^{1/2}\
\prod _{k=0}^{\frac {n-3}2} (-\Delta + k(n-k-1)).
\end{cases}
\end{equation}
Using the method of moving planes, it is shown in \cite{CY-6} that
all solutions of the (pseudo-) differential equation:
\begin{equation}
\mathbb P_n w + (n-1)!=(n-1)!e^{nw}
\end{equation}
are given by actions of the conformal group of $S^n$.
As a consequence, we derive (\cite{CY-5}) the sharp version of a Moser-Trudinger
inequality for spheres in general dimensions.
This inequality is equivalent to Beckner's inequality (\cite{Be}).
\begin{equation}
\log \frac{1}{|S^n|}\int _{S^n}e^{nw} dv \leq \frac{1}{|S^n|}\int _{S^n}
(nw + \frac{n}{2 (n-1 )!} w \mathbb P_n (w)) dv,
\end{equation}
and equality holds if and only if $e^{nw}$ represents the Jacobian of
a conformal transformation of $S^n$.

In a recent preprint, S. Brendle is able to derive a general existence result
for the prescribed $Q$-curvature equation under natural conditions:

\begin{Thm}{\rm (\cite{B})}
For a compact manifold $(M^{2m},g)$ satisfying
\item{(i)} $P_{2m}$ be positive except on constants,
\item{(ii)} $\int_M Q_g dv_g < C_{2m}$
where $C_{2m}$ represents the value of the corresponding
$Q$-curvature integral on the standard sphere $(S^{2m}, g_c)$,
the equation $P_{2m}w+ Q \, = \, Q_w e^{2mw}$ has a solution with
$ Q_w $ given by a constant.
\end{Thm}

Brendle's remarkable argument uses a $2m$-th order heat flow method in which again inequality
of Adams (\cite{A}) (the only available tool) is used.

In another recent development, the $n$-th order $Q$-curvature integral
can be interpreted as a renormalized volume of the conformally compact
manifold $(N^{n+1},h)$ of which $(M^n,g)$ is the conformal infinity.
In particular, Graham-Zworski (\cite{GZ}) and Fefferman-Graham (\cite{FG-2})
have given in the case $n$ is an even integer, a spectral theory interpretation
to the $n$-th order $Q$-curvature integral that is intrinsic to the boundary
conformal structure. In the case $n$ is odd, such an interpretation is still
available, however it may depend on the conformal compactification.
%yang June 13, 02

\section{\hskip -1em . \hskip 0.3em Boundary operator, Cohn-Vossen inequality}
\setzero \vskip-5mm \hspace{5mm} \setcounter{Thm}{0}

To develop the analysis of the $Q$-curvature equation, it is
helpful to consider the associated boundary value problems.
In the case of compact surface with boundary $(N^2, M^1, g)$ where
the metric $g$ is defined on $N^2 \cup M^1$; the Gauss-Bonnet
formula becomes
\begin{equation}
2 \pi \chi(N) = \int_N K \,\,  dv + \oint_M k \,\,  d\sigma,
\end{equation}
where $k$ is the geodesic curvature on $M$. Under conformal
change of metric $g_w$ on $N$, the geodesic curvature changes according to the equation
\begin{equation}
\frac{\partial}{\partial n} w + k \, = \, k_w e^{w} \,\,\, \text {on M}.
\end{equation}
Ray-Singer-Polyakov log-determinant formula
has been generalized to compact surface with boundary
and the extremal metric of the formula has been studied by Osgood-Phillips-Sarnak
(\cite {OPS-2}). The role played by the Onofri inequality is the
classical Milin-Lebedev inequality:
\begin{equation}
\log \oint_{S^1} e ^{(w - \bar w)} \frac {d \theta}{ 2 \pi} \leq
\frac{1}{4} \left( \int_D w (- \, \Delta w) \frac {dx}{\pi} \, +
\, 2
 \oint _{S^1}
 w \frac{\partial w }{\partial n} \frac {d \theta}{ 2 \pi} \right),
\end{equation}
where $D$ is the unit disc on $\mathbb R^2$ with the flat metric $dx$, and $n$ is the unit outward normal.

One can generalize above results to
four manifold with boundary $(N^4, M^3 ,g)$; with the role
played by $ (- \Delta,  \frac{\partial}{\partial n})$ replaced by
$(P_4, P_3)$ and with  $(K, k)$ replaced by
$(Q, T)$; where $P_4$ is the Paneitz opertor and $Q$ the
curvature discussed in section 2; and where $P_3$ is the boundary operator
constructed by Chang-Qing (\cite {CQ-1}).  The key property of $P_3$ is
that
it is conformally covariant
 of bidegree $(0, 3)$,
when operating on functions defined on the boundary of
 compact $4$-manifolds; and under conformal change of metric
 $\bar g= e^{2w}g$ on $N^4$
we have at the boundary $M^3$
\begin{equation}
P_3 w + T \, = \, T_w e^{3w}.
\end{equation}
We refer the reader to \cite{CQ-1} for the precise definitions of
$P_3$ and $T$ and will here only mention that on $(B^4, S^3, dx)$,
where $B^4$ is the unit ball in $\mathbb R^4$, we have
\begin{equation}
P_4 = (- \Delta )^2, \,\, P_3 = -\left( \frac{1}{2}\,\,  \frac{\partial}{\partial n} \,\, \Delta + \tilde \Delta \frac{\partial}{\partial n} + \tilde \Delta \right) \,\,\, \,\,\,
\text {and} \,\,\,
T = 2,
\end{equation}
where $\tilde \Delta$ is the intrinsic boundary Laplacian on $M$.

In this case the Gauss-Bonnet-Chern formula may be expressed as:
\begin{equation}
4 \pi ^2 \chi (N)= \int _N  (\, Q \, + \frac{1}{8}|W|^2 )\,\, dv
\,\, + \oint _M (T + {\cal L}) \,\,  d\sigma,
\end{equation}
where $\cal L$ is a third order boundary curvature invariant that
t ransforms by scaling under conformal change of metric. The
analogue of the sharp form of the Moser-Trudinger inequality for
the pair $(B^4,S^3,\, dx)$ is given by the following analogue of
the Milin-Lebedev inequality:

\begin{Thm}{\rm (\cite{CQ-2})}
Suppose $w \in C^\infty (\bar B^4)$.
Then
\begin{eqnarray}
 & & \log \left\{ \frac{1}{2\pi^2} \oint_{S^3} e^{3(w-\bar w)} d \sigma
 \right\} \nonumber \\
  & \leq&
\frac{3}{16\pi^2} \left\{ \int_{B^4} w \Delta^2 w dx \,+ \,
 \oint_{S^3} \left( 2 w P_3 w - \frac {\partial w}
{\partial n} + \frac { \partial ^2 w}{\partial n^2} \right) \, d \sigma
 \right\},
\end{eqnarray}
under the boundary assumptions
$ \frac{\partial w}{\partial n} |_{S^3} = e^{w} -1 $ and $\int_{S^3} R_w d\sigma_{g_w}  = \int_{S^3} R d\sigma $
where $R$ is the scalar curvature of $S^3$.
Moreover the equality holds if and only if $e^{2w} dx$ on
 $B^4$ is
isometric to the standard metric via a conformal transformation of the pair
$(B^4,S^3, dx)$.
\end{Thm}

The boundary version (3.6) of the Gauss-Bonnet-Chern formula can be
used to give an extension of the well known Cohn-Vossen-Huber formula.
Let us recall (\cite{CV}, \cite{Hu})
that a complete surface $(N^2,g)$ with Gauss curvature in $L^1$
 has a conformal compactification $\bar {N}= N \cup \{q_1, ... ,q_l\}$
as a compact Riemann surface and
\begin{equation}
2 \pi \chi (N)= \int _{N} KdA + \sum _{k=1}^{l} \nu_k,
\end{equation}
where at each end $q_k$, take a conformal coordinate disk $\{|z| < r_0\}$
with $q_k$ at its center, then $\nu_k$ represents the following limiting
isoperimetric constant:
\begin{equation}
\nu _k = \lim _{r \rightarrow 0} \frac {{Length(\{|z|=r\})}^2} {2 Area(\{r<|z|<r_0 \})}.
\end{equation}

This result can be generalized to dimension $n=4$ for locally conformally flat
metrics. In general dimensions, Schoen-Yau (\cite{SY}) proved that locally
conformally flat
metrics in the non-negative Yamabe class has injective development map into
the standard spheres as domains whose complement have small Hausdorff dimension
(at most $\frac{n-2}{2}$).
It is possible to further constraint the topology as well as
the end structure of such manifolds by imposing the natural condition that the $Q$-curvature be in $L^1$.

\begin{Thm}{\rm (\cite{CQY-1}, \cite{CQY-2})}
Suppose $(M^4,g)$ is a complete conformally flat manifold, satisfying the
conditions:
\newline
(i) The scalar curvature $R_g$ is bounded between two positive constants and $\nabla_g R_g$
is also bounded;
\newline
(ii) The Ricci curvature is bounded below;
\newline
(iii) $\int_M |Q_g|dv_g < \infty$;
\newline
then
\newline
(a) if $M$ is simply connected, it is conformally equivalent to
$S^4-\{q_1, ... ,q_l\}$
and we have
\begin{equation}
 4 \pi^2 \,\, \chi(M) = \int_M Q_g \,\, dv_g\,\, +\,\,  4 \pi^2 l\,\,\,  ;
\end{equation}
\newline
(b) if $M$ is not simply connected, and we assume in addition that its fundamental
group is realized as a geometrically finite Kleinian group, then we conclude
that $M$ has a conformal compactification $\bar M= M \cup \{q_1, ... ,q_l\}$
and  equation (3.10) holds.
\end{Thm}

This result gives a geometric interpretation to the $Q$-curvature integral
as measuring an isoperimetric constant.
There are two elements in this argument. The first
 is to view the $Q$-curvature integral over
sub-level sets of the conformal factors as the second  derivative
with respect to $w$ of the corresponding
volume integral. This comparison is made possible by making use of the
formula (3.4). A second element is
an estimate for conformal metrics $e^{2w}|dx|^2$
defined over domains  $\Omega \subset \Bbb R ^4$ satisfying the
conditions of Theorem 3.2
must have a uniform blowup rate near the boundary:
\begin{equation}
e^{w(x)} \cong \frac{1}{d(x, \partial \Omega)}.
\end{equation}
 This result has an appropriate generalization to higher
even dimensional situation, in which one has to impose additional curvature
bounds to control the lower order terms in the integral. One such an extension
is obtained in the thesis of H. Fang (\cite{F}).

It remains an interesting question how to extend this analysis to include
the case when the dimension is an odd integer.

\section{\hskip -1em . \hskip 0.3em Fully nonlinear equations in
conformal geometry in dimension four} \setcounter{Prop}{0}
\setcounter{Thm}{0} \setzero \vskip-5mm \hspace{5mm}

In dimensions greater than two, the natural curvature invariants in
conformal geometry
are the Weyl tensor $W$, and the Weyl-Schouten tensor $A=Ric - \frac{R}{2(n-1)}g$
that occur in the
 decomposition of the curvature tensor;
where $Ric$ denotes the Ricci curvature tensor:
\begin{equation}
Rm=W \oplus  \frac{1}{n-2} A \bcw  g.
\end{equation}
Since the Weyl tensor $W$ transforms by scaling under conformal change
$g_w= e^{2w}g$,
only the Weyl-Schouten tensor depends on the derivatives of the
conformal factor. It is thus natural to consider $\sigma_k(A_g)$
the k-th symmetric function of the eigenvalues
of the Weyl-Schouten tensor $A_g$ as curvature invariants of
the conformal metrics.
As a differential invariant of the conformal factor $w$, $\sigma_k(A_{g_w})$ is a
fully nonlinear expression involving the Hessian and the gradient of
the conformal factor $w$. We have abbreviating $A_w$ for $A_{g_w}$:
\begin{equation}
A_w= (n-2) \{- \nabla ^2w +  dw\otimes dw- \frac{|\nabla w|^2}{2} \}
+ A_g.
\end{equation}
The equation
\begin{equation}
\sigma _k(A_w)= 1
\end{equation}
is a fully nonlinear version of the Yamabe equation. For example, when
$k =1$, $\sigma_1(A_g) = \frac{n-2}{ 2(n-1)} R_{g}$, where $ R_g$ is the
scalar curvature of $(M, g)$ and  equation (4.3) is the Yamabe equation
which we have discussed in section 1. When $ k =2 $, $\sigma_2 (A_g) =
\frac {1}{2} (|Trace \,\, A_g|^2 - |A_g|^2) = \frac{n}{8 (n-1)} R^2 -
\frac{1}{2} |Ric|^2 $. In the case when $k =n$, $\sigma_n (A_g) =
determinant \,\, of \,\, A_g$, an equation of Monge-Ampere type. To
illustrate that (4.3) is a fully non-linear elliptic equation, we have
for example when $n=4$,
\begin{equation}
 \begin{aligned}
\sigma_2(A_{g_w}) e^{4w} \, =&\, \sigma_2 (A_g) \, + 2 ((\Delta w)^2
- \, |\nabla ^2 w|^2 \\
+& \,  (\nabla w,\nabla |\nabla w|^2)
+  \Delta w |\nabla w |^2 )\,\, \\
+& \text {lower order terms},
\end{aligned}
\end{equation}
where all derivative are taken with respect to the $g$ metric.

For a symmetric $n\times n$ matrix $M$, we say $M \in \Gamma _k^+$
in the sense of Garding  (\cite {Ga})
if
$\sigma _k(M) >0$ and $M$ may be joined to the identity matrix by a path
consisting entirely of matrices $M_t$ such that $\sigma _k(M_t) >0$.
There is a rich literature conerning the equation
\begin{equation}
\sigma _k(\nabla ^2 u)= \, f ,
\end{equation}
for a positive function $f$. In the case when $M = ( \nabla^2 u )$ for
convex functions $u$ defined on the Euclidean domains, regularity theory
for equations of $\sigma_k(M)$ has been well established for $M \in
\Gamma _k^+$ for Dirichlet boundary value problems by
Caffarelli-Nirenberg-Spruck (\cite{CNS-2}); for a more general class of
fully non-linear elliptic equations not necessarily of divergence form
by Krylov (\cite {Kr}), Evans (\cite {E-1}) and for Monge-Ampere
equations by Pogorelov (\cite{Pog}) and by Caffarelli (\cite{Ca-1}). The
Monge-Ampere equation for prescribing the Gauss-Kronecker curvature for
convex hypersurfaces has been studied by Guan-Spruck (\cite{GS}).
 Some of the techniques in these work
can be modified to study equation (4.3) on manifolds. However there are
features of the equation (4.3) that are distinct from the equation
(4.5). For example, the conformal invariance of the equation (4.3)
introduces a non-compactness due to the action of the conformal group
that is absent for the equation (4.5).

When $ k \neq \frac{n}{2} $ and the manifold $(M,g)$ is locally
conformally flat, Viaclovsky (\cite {V-1}) showed that the equation
(4.3) is the Euler equation of the variational functional $\int
\sigma_k(A_{g_w})dv_{g_w}$. In the exceptional case $k = n/2$, the
integral $\int \sigma _k(A_{g})dv_{g}$ is a conformal invariant. We say
$g \in \Gamma _k^+$ if the corresponding Weyl-Schouten tensor $A_g(x)
\in \Gamma _k^+$ for every point $x \in M$.  For $k=1$ the Yamabe
equation (1.10) for prescribing scalar curvature is a semilinear one;
hence the condition for $g \in \Gamma _1^+$ is the same as requiring the
operator $L_g = -\frac{4(n-1)}{n-2}\Delta_g  + R_g $ be a positive
operator. The existence of a metric with $g \in \Gamma _k^+$ implies a
sign for the curvature functions (\cite {GV}, \cite{CGY-3}, \cite
{GVW}).
\begin{Prop}
On $(M^n, g)$,
\newline
(i) When $n=3$ and $\sigma_2 (A_g)>0 $, then either $R_g >0$ and
the sectional curvature of g is positive or $R_g < 0$ and
the sectional curvature of g is negative on $M$.
\newline
(ii) When $n=4$ and $\sigma_2 (A_g)>0 $, then either $R_g >0$ and
$Ric_g >0$ on $M$ or $R_g <0$ and $Ric_g < 0$ on $M$.
\newline
(iii) For general $n$ and $A_g \in \Gamma _k^+$ for some
$k \geq \frac{n}{2}$, then $Ric_g > 0$.
\end{Prop}

In dimension 3, one can capture all metrics with constant
sectional curvature (i.e. space forms) through the study of $\sigma_2$.

\begin{Thme} {\rm (\cite{GV})}{\it
\,\, On a compact 3-manifold, for any Riemannian metric $g$,\, denote
${\cal F}_2 [g] = \int_M \sigma_2 (A_g) dv_g$.  Then a metric $g$ with $
{\cal F}_2 [g] \geq 0$ is critical for the functional ${\cal F}_2 $
restricted to class of metrics with volume one if and only if $g$ has
constant sectional curvature.}
\end{Thme}

The criteria for existence of a conformal metric $g \in \Gamma _k^+$ is not
as easy for $k>1$ since the equation is a fully nonlinear one. However
when $n=4, k=2$ the invariance of the integral $\int \sigma _2(A_g)dv_g$
is a reflection of the Chern-Gauss-Bonnet formula
\begin{equation}
8\pi^2 \chi (M)= \int _M(\sigma_2 (A_g) + \frac{1}{4}|W_g|^2)dv_g.
\end{equation}

In this case it is possible to find a criteria:
\begin{Thme}{\rm (\cite {CGY-3})}{\it
For a closed 4-manifold $(M, g)$ satisfying the following conformally invariant
conditions:
\newline\noindent
(i) $Y(M, g)\, > \, 0, $ and
\newline\noindent
(ii) $\int \sigma _2(A_g)dv_g >0$;\newline\noindent then there exists a
conformal metric $g_w \in \Gamma _2^+$.}
\end{Thme}

\noindent {\bf Remark.} In dimension four, the condition $g \in
\Gamma _2^+$ implies that $R>0$ and Ricci is positive everywhere.
Thus such manifolds have finite fundamental group. In addition,
the Chern-Gauss-Bonnet formula and the signature formula shows
that this class of 4-manifolds satisfy the same conditions as that
of an Einstein manifold with positive scalar curvatures. Thus it
is the natural class of 4-manifolds in which to seek an Einstein
metric.

The existence result depends on the solution of a family of fourth order equations
involving the Paneitz operator (\cite{Pa}), which we have discussed in section 2. In the following we briefly outline
this connection. Recall that in  dimension four, the Paneitz operator $P$
a fourth order curvature called the Q-curvature:
\begin{equation}
P_g w + 2Q_g \, = \, 2 Q_{g_w} e^{4w}.
\end{equation}
The relation between $Q$ and $\sigma_2(A)$ in dimension 4 is given
by
\begin{equation}
Q_g = \frac{-1}{12}\Delta R_g + \frac{1}{2}\sigma_2 (A_g).
\end{equation}

In view of the existence results of Theorem 2.1 and Theorem 2.2, it is
natural to find a solution of
\begin{equation}
\sigma_2 (A_g) \, = \, f
\end{equation}
for some positive function $f$. It turns out that it is natural to
choose $ f =  c |W_g|^2$ for some constant $c$ and to use the continuity
method to solve the family of equations
\begin{equation}
(*)_\delta : \, \,  \,\,\,\, \,\,\,\,\,\,\,\,  \sigma_2(A_g)= \frac{\delta}{4} \Delta_g R_g - 2 \gamma |W_g|^2
\end{equation}
where $\gamma$ is chosen so that $\int \sigma _2(A_g)dv_g=- 2
\gamma \int|W_g|^2 dv_g$, for $\delta \in (0, 1]$ and let $\delta$
tend to zero.

Indeed when $\delta =1 $, solution of (4.10) is a special case of
an extremal metric of the log-determinant type functional $F[w]$
in Theorem 2.2, where we choose $\gamma_2 = 1 $, $\gamma_3 = \frac
{1}{24}$, we then choose $\gamma= \gamma_1$ so that $\bar k = 0$.
Notice that in this case, the assumption (ii) in the statement of
Theorem 4.3 implies that $\gamma < 0 $. When $\delta = \frac
{2}{3} $, equation (4.10) amounts to solving the equation
\begin{equation}
Q_g \, = \,  - \gamma |W_g|^2,
\end{equation}
which we can solve by applying Theorem 2.1. Thus the bulk of the analysis consist in
obtaining apriori estimates of the solution as $\delta$
tends to zero, showing essentially that in the equation the term
$\frac{\delta}{4}\Delta R$ is small in the weak sense.
The proof ends by first modifying the function $|W|^2$ to make it strictly
positive and by then applying the Yamabe flow  to the metrics
$g_{\delta}$
 to show that for sufficiently small $\delta$ the smoothing
provided by the Yamabe flow yields a metric $g \in \Gamma _2^+$.

The equation (4.3) becomes meaningful for 4-manifolds which admits a metric
$g \in \Gamma _2^+$. In the article (\cite{CGY-4}), when the manifold
$(M, g)$ is not conformally equivalent to $(S^4, g_c)$, we provide apriori estimates
for solutions of the equation (4.9)
where $f$ is a given positive smooth function. Then we apply the degree theory for fully non-linear elliptic equation
% of Y. Li {\cite {L-1} )
to the following
1-parameter family of equations
\begin{equation}
\sigma _2(A_{g_t})= tf +(1-t)
\end{equation}
to deform the original metric to one with constant $\sigma _2(A_g)$.

In terms of geometric application, this circle of ideas may be applied to
characterize a number of interesting conformal classes in terms of the
the relative size of the conformal invariant $\int \sigma _2(A_g)dV_g$
compared with the Euler number.

\begin{Thme}
{\rm (\cite{CGY-6})}{\it \,\, Suppose $(M,g)$ is a closed 4-manifold
with $Y(M, g) \, > \, 0$.

(I) If $\int _M \sigma _2(A_g)dv_g > \frac {1}{4} \int _M |W_g|^2 \,\,
dv_g$, then $M$ is diffeomorphic to $(S^4, g_c)$ or $(\mathbb RP^4,
g_c)$.

(II) If $M$ is not diffeomorphic to $(S^4, g_c)$ or $(\mathbb RP^4,
g_c)$ and $\int _M \sigma _2(A_g)dv_g = \frac {1}{4} \int _M |W_g|^2
\,\, dv_g$, then either

(a) $(M, g)$ is conformally equivalent to $(\mathbb CP^2, g_{FS})$, or

(b) $ (M, g)$ is conformal equivalent to $( (S^3 \times S^1)/\Gamma,
g_{prod}). $}
\end{Thme}
\noindent {\bf Remark.} The theorem above is an $L^2$ version of an
earlier result of Margerin \cite{Ma}. The first part of the theorem
should be compared to a result of Hamilton (\cite {H-1}); where he
pioneered the method of Ricci flow and established the diffeomorphism of
$M^4$ to the 4-sphere under the assumption that the curvature operator
be positive.

This first part of Theorem 4.4 applies the existence argument to find
a conformal metric $g'$ which satisfies the pointwise inequality

\begin{equation}
\sigma_2 (A_{g'}) > \frac {1}{4}  |W_{g'}|^2.
\end{equation}
The diffeomorphism assertion follows from Margerin's (\cite{Ma}) precise
convergence result for the Ricci flow: such a metric will evolve under
the Ricci flow to one with constant curvature. Therefore such a
manifold is diffeomorphic to a quotient of the standard $4$-sphere.

For the second part
of the assertion, we argue that if such a manifold is not
diffeomorphic to the 4-sphere, then the conformal structure realizes
the minimum of the quantity $\int |W_g|^2 dv_g$, and hence its Bach tensor
vanishes. There are two possibilities depending on whether the
Euler number is zero or not. In the first case, an earlier result
of Gursky (\cite{G-1}) shows the metric is conformal to that of the space
 $S^1 \times S^3$. In the second case, we solve the equation
\begin{equation}
\sigma _2(A_{g'}) = \frac {1 - \epsilon }{4}  |W_{g'}|^2 + \, C_{\epsilon},
\end{equation}
where $C_{\epsilon}$ is a constant which tend to zero as $\epsilon$ tend
to zero. We then let $\epsilon$ tends to zero. We obtain in the limit a
$C^{1,1}$ metric which satisfies the equation on the open set $\Omega =
\{x| W(x) \neq 0\}$:
\begin{equation}
\sigma _2(A_{g'}) = \frac {1}{4}  |W_{g'}|^2.
\end{equation}
Then a Lagrange multiplier computation
shows that the curvature
tensor of the limit metric agrees with that of the Fubini-Study metric
on the open set where $W\neq 0$. Therefore $|W_{g'}|$ is a constant on $\Omega$
thus $W$ cannot vanish at all. It follows from the Cartan-Kahler
theory that
the limit metric agrees with the Fubini-Study metric of $\mathbb CP^2$
everywhere.

There is a very recent work of A. Li and Y. Li (\cite{LL}) extending
work of (\cite{CGY-5}) to classify the entire solutions of the equation
$\sigma _k(A_g)=1$ on $\mathbb R ^n$ thus providing apriori estimates
for this equation in the locally conformally flat case. There is also a
very recent work (\cite{GW}) on the heat flow of this equation, we have
(\cite{CY-8}) used this flow to derive the sharp version of the
Moser-Onofri inequality for the $\sigma _{\frac{n}{2}}$ energy for all
even dimensional spheres. In general, the geometric implications of the
study of $\sigma_{k}$ for manifolds of dimension greater than four
remains open.

\newcommand{\namelistlabel}[1] {\mbox{#1}\hfil}
\newenvironment{namelist}[1]{%
\begin{list}{}
{\let\makelabel\namelistlabel
\settowidth{\labelwidth}{#1}
\setlength{\leftmargin}{1.1\labelwidth}}
}{%
\end{list}}

%\begin{document}

%\newpage

\label{lastpage}

\end{document}